\documentclass{article}
\usepackage{amssymb}
\usepackage{amsfonts}
\usepackage{amsmath}
\usepackage{latexsym}

\newtheorem{theorem}{Theorem}[section]
\newtheorem{corollary}[theorem]{Corollary}
\newtheorem{lemma}[theorem]{Lemma}
\newtheorem{proposition}[theorem]{Proposition}
\newtheorem{example}[theorem]{Example}
\newtheorem{problem}[theorem]{Problem}
\newtheorem{remark}[theorem]{Remark}
\newtheorem{definition}[theorem]{Definition}

\newcommand{\FF}{\mathbb F}
\newcommand{\CC}{\mathbb C}

\def\Gl{\mathop{\rm Gl}\nolimits}
\def\deg{\mathop{\rm deg }\nolimits}
\def\rank{\mathop{\rm rank}\nolimits}

\newcommand{\se}{\ensuremath{\stackrel{s.e.}{\sim}}}
\def\lcm{\mathop{\rm lcm }\nolimits}

\newcommand{\ba}{\mathbf a}
\newcommand{\bd}{\mathbf d}
\newcommand{\bg}{\mathbf g}
\newcommand{\bc}{\mathbf c}
\newcommand{\bu}{\mathbf u}
\newcommand{\bv}{\mathbf v}
\newcommand{\bw}{\mathbf w}
\newcommand{\be}{\mathbf e}
\newcommand{\bb}{\mathbf b}

\newcommand{\MSC}[1]{\textbf{\textit{MSC}} #1}
\title{Rank-one perturbations of matrix pencils}
\author{Itziar Baraga\~na\footnote{Departamento de Ciencia de la Computaci\'on e I.A.,
Facultad de Inform\'atica, Universidad del Pa\'{\i}s Vasco, UPV/EHU, 
 Donostia-San Sebasti\'an, Spain, e-mail: itziar.baragana@ehu.eus. 
 Partially supported by  
``Ministerio de Econom\'{\i}­a, Industria y Competitividad (MINECO)'' of Spain and ``Fondo Europeo de Desarrollo Regional (FEDER)'' of EU through grants MTM2017-83624-P and MTM2017-90682-REDT, and by UPV/EHU through grant GIU16/42.},
Alicia Roca\footnote{
Departamento de  Matem\'atica Aplicada,
 IMM, Universitat Polit\`ecnica de Val\`encia,   Valencia, Spain, e-mail:
aroca@mat.upv.es.
 Partially supported by  ``Ministerio de Econom\'{\i}­a, Industria y Competitividad (MINECO)'' of Spain and ``Fondo Europeo de Desarrollo Regional (FEDER)'' of EU through grants MTM2017-83624-P and MTM2017-90682-REDT.
}
}
\date{}

\begin{document}

\maketitle

\begin{abstract}
We solve the problem of characterizing the Kronecker structure of a  matrix pencil obtained by  a rank-one  perturbation of another matrix pencil. The results hold over  arbitrary fields.
\end{abstract}


\MSC
  15A21, 15A22, 47A55

\section{Introduction}

Given a matrix pencil $A(s)=A_0+sA_1 \in \FF^{n\times m}$, the  rank perturbation  problem consists in characterizing the Kronecker structure of  $A(s)+P(s)$, where  $P(s)$ is a matrix pencil of  bounded rank.

The Kronecker structure of a matrix pencil is determined by the complete system of invariants for the strict equivalence of matrix pencils, i.e., the invariant factors, infinite elementary divisors, and row and column minimal indices. For regular matrix pencils the Kronecker structure is known as the  Weierstrass structure and is determined only by  the invariant factors and the infinite elementary divisors.
In particular, the Jordan structure of a square matrix is defined  by the Weierstrass structure of the associated characteristic pencil, which is a regular pencil without infinite elementary divisors.
Analogously, the feedback invariants of a pair of matrices, i.e., the invariant factors and the column (or row) minimal indices  are the Kronecker invariants of the associated characteristic pencil.  

In the last decades rank perturbations of matrix pencils have been analyzed   in many papers  from different  approaches. The problem has been studied  generically, i.e.,  when the perturbation $P(s)$ belongs to an open and dense subset of the set of pencils of rank less than or equal to $r$, for a given integer $r$. In other cases, the pencil $P(s)$ is an arbitrary perturbation belonging to the whole set of pencils of  rank less than or equal to $r$.  In this paper we follow the second approach. 

From a generic point of view, changes in the Jordan structure of a square constant matrix or in the Weierstrass structure of a regular pencil corresponding to a fixed eigenvalue after low rank perturbations have been studied, among others,   in \cite{Batzke14, BaMeRaRo16, TeDo16, TeDoMo08, MeMeRaRo11, MoDo03, Sa02, Sa04}.
See also the references therein. 
 
The case where the perturbation is an arbitrary square matrix $P$ or a regular matrix pencil $P(s)$ has also been studied by several authors. 
For square constant matrices and a constant perturbation of bounded rank $r$, a  solution is given in \cite{Silva88_1} and \cite{Za91}. For $r=1$ the problem was already solved in \cite{Th80}.
The case where the perturbation has fixed rank  has been solved in   \cite{Silva88_1}   over algebraically closed fields. 

For regular pencils  the problem has been studied  for $r=1$ in \cite{GeTr17}. For arbitrary perturbations of bounded rank the problem has been solved in  \cite{BaRo18}, and for perturbations of  fixed rank in \cite{BaRo19}. In both cases the solutions obtained do not involve  any condition on the rank of the type ``low-rank'', and the results hold for  fields having sufficient number of elements (fields  requiring just the condition that at least one element of the field or the point at infinity is not included neither in the spectrum of the original pencil nor in the perturbed one).

There is less literature dealing with the case of singular pencils. The problem is more difficult, since the  row and column minimal indices of the pencils are involved.
For  non full rank pencils the change of the four types of invariants under generic low rank perturbations has been characterized in  \cite{TeDo07}.
For square singular pencils, in \cite{LeMaPhTrWi18} the authors represent pencils via linear relations and obtain bounds for the number of Jordan  chains which may change under an arbitrary rank-one perturbation.
The problem of characterizing the feedback equivalence  invariants of a pair of matrices, i.e., the Kronecker invariants of the associated characteristic pencil, under  a constant perturbation of bounded rank  is solved in \cite{DoSt14}. Here, the authors find the solution relating   the problem  to a matrix pencil completion problem.

\medskip
 
In this paper we study arbitrary rank-one  perturbations of matrix pencils. We solve the problem transforming it into a matrix pencil completion problem. The solution obtained holds for arbitrary fields.

\medskip

The paper is organized as follows. In Section \ref{secpreliminaries} we introduce the notation, basic definitions and preliminary results. 
In Section \ref{secproblem} we establish the problem which we are going to study and relate it to a matrix pencil completion problem.
Then, in Section \ref{seccompletion} we introduce previous  results  about completion of matrix pencils which will be needed later. In Section \ref{sectechnical}
some thecnical lemmas are proved. In Section \ref{secmain},  a solution to the stated rank-one perturbation problem is given in Theorem \ref{maintheogen}.
Finally, in Section \ref{secconclusions} we summarize the main contributions of the paper.

\section{Preliminaries}
\label{secpreliminaries}

Let $\FF$ be a field. $\FF[s]$ denotes the ring of polynomials in the indeterminate $s$ with coefficients in $\FF$ and  $\FF[s, t]$  the ring of polynomials in two
variables $s, t$ with coefficients in $\FF$.
We denote by $\FF^{p\times q}$,  $\FF[s]^{p\times q}$ and $\FF[s, t]^{p\times q}$ the vector spaces  of $p\times q$ matrices with elements in $\FF$, $\FF[s]$ and $\FF[s, t]$, respectively.
$\Gl_p(\FF)$ will be the general linear group of invertible matrices
in $\FF^{p \times p}$. 

\medskip

Given a polynomial matrix $G(s)\in \FF[s]^{p\times q}$, the {\em degree} of $G(s)$, denoted by $\deg(G(s))$, is the  maximum  of  the degrees of its entries. The {\em normal rank} of $G(s)$, denoted by $\rank (G(s))$,  is the order of the largest nonidentically zero minor of $G(s)$, i.e., it is the rank of $G(s)$ considered as a matrix on the field of fractions of $\FF[s]$. 

\medskip

A {\em matrix pencil} is a  polynomial matrix $G(s)\in \FF[s]^{p\times q}$
such that $\deg(G(s))\leq1$. The pencil is  {\em regular} if $p=q$ and $\det(G(s))$ is not the zero polynomial. Otherwise it is  {\em singular}.

Two matrix pencils $G(s)=G_0+sG_1, H(s)=H_0+sH_1\in \FF[s]^{p\times q}$ are {\em strictly equivalent} ($G(s)\se H(s)$) 
if there exist invertible matrices $Q\in \Gl_p(\FF)$,   $R\in \Gl_q(\FF)$ such that 
$G(s)=QH(s)R$. 

Given the pencil $G(s)=G_0+sG_1 \in \FF[s]^{p\times q}$ of $\rank G(s)=n$, a complete system of invariants for  the strict equivalence  of matrix pencils is formed by a chain of homogeneous polynomials $\Gamma_1(s,t) \mid \dots \mid \Gamma_n(s,t),\ \Gamma_i(s,t) \in \FF[s,t], \ 1\leq i\leq n$, called the  {\em homogeneous invariant  factors}, and two collections of nonnegative integers $c_1\geq \dots \geq c_{q-n}$ and $u_1\geq \dots \geq u_{p-n}$, called  the {\em column and row minimal indices} of the pencil, respectively. In turn, the homogeneous invariant factors are determined by a chain of polynomials  $\gamma_1(s)\mid \ldots \mid \gamma_n(s)$ in $\FF[s]$, called the invariant factors, and a chain of polynomials  $t^{k_1}\mid \ldots \mid t^{k_n}$ in $\FF[t]$, called  the infinite elementary divisors. 
In fact, we can write 
$$\Gamma_i(s,t)=t^{k_i}t^{\deg(\gamma_i(s))}\gamma_i\left(\frac st\right), \ 1\leq i\leq n.$$ 
The associated canonical form is the Kronecker canonical form. For details see  \cite[Ch. 2]{Friedland80} or \cite[Ch. 12]{Ga74} for infinite fields, and  \cite[Ch. 2]{Ro03} for arbitrary fields.
In what follows we will work with the homogeneous invariant factors.
We will take $\Gamma_i(s,t)=1$ ($\gamma_i(s)=1$) whenever $i<1$ and $\Gamma_i(s,t)=0$ ($\gamma_i(s)=0$) when $i>n$.
The sum of the degrees of the homogeneous invariant factors  plus the sum of the minimal indices is equal to the rank of the pencil. 
Also, if $T(s)=G(s)^T$, then $G(s)$ and $T(s)$ share the homogeneous invariant factors and have interchanged minimal indices, i.e., the column (row) minimal indices of $T(s)$ are the row (column) minimal indices of $G(s)$.

Observe  that if $G(s)\in \FF[s]^{p\times q}$ and $\rank (G(s))=p$ ($\rank (G(s))=q$), then $G(s)$ does  not have row (column) minimal indices. As a consequence, 
the invariants for the strict equivalence of  regular matrix pencils are reduced to  the homogeneous invariant factors.

\medskip

In this paper we study the Kronecker structure of arbitrary pencils perturbed by pencils of rank one.
A matrix pencil of rank one allows a very simple decomposition (see \cite{GeTr17} for $\FF=\CC$). In the next proposition we analyze this decomposition  for arbitrary fields, depending on the Kronecker structure of the pencil.

\begin{proposition}\label{pencil1}
Let $P(s)\in\FF[s]^{p\times q}$ be a matrix pencil of $\rank P(s)=1$. 
\begin{enumerate}
\item If  $P(s)$ has a nontrivial invariant factor, then there exist nonzero vectors 
$u\in \FF^{p}$, $\bar v\in \FF^q$ and nonzero pencils $\bar u(s)\in \FF[s]^{p}$, $v(s)\in \FF[s]^{q}$
such that 
$$
P(s)=uv(s)^T=\bar u(s)\bar v^T.
$$
\item If $P(s)$ has an infinite elementary divisor, then there exist nonzero vectors $u \in \FF^{p}$, $v\in \FF^q$
such that
$$
P(s)=uv^T.
$$
\item If $P(s)$ has a positive column minimal index, then there exist a nonzero vector
$u\in \FF^{p}$ and a nonzero pencil  $v(s)\in \FF[s]^{q}$
such that
$$
P(s)=uv(s)^T.
$$
\item If $P(s)$ has a positive row minimal index, then there exist a nonzero vector $v\in \FF^{q}$ and a nonzero pencil  $u(s)\in \FF[s]^{p}$
such that 
$$
P(s)=u(s)v^T.
$$
\end{enumerate}
\end{proposition}

{\it Proof.}
Let $P_c(s)$ be the Kronecker canonical form of $P(s)$. 
Then, there exist
$Q=\begin{bmatrix}q_1&\dots&q_p\end{bmatrix}\in \Gl_{p}(\FF)$ and 
$R=\begin{bmatrix}r_1^T\\\vdots\\r_q^T\end{bmatrix}\in \Gl_{q}(\FF)$
such that
$P(s)=QP_c(s)R$.
\begin{enumerate}
\item If $P(s)$ has a nontrivial invariant factor $s+\lambda$, $\lambda \in \FF$, then 
$P_c(s)=\begin{bmatrix}
s+\lambda &0 \\0&0
\end{bmatrix}\in \FF[s]^{p\times q}$.
Hence,
$P(s)=((s+\lambda)q_1)r_1^T=q_1((s+\lambda)r_1)^T.$
\item If $P(s)$ has an infinite elementary  divisor, then 
$P_c(s)=\begin{bmatrix}
1&0\\
0&0\end{bmatrix}\in \FF[s]^{p\times q}.$
Therefore,
$P(s)=q_1r_1^T$.
\item If $P(s)$ has a positive column minimal index, then
$P_c(s)=
\begin{bmatrix}
s&1 &0\\
0&0 &0\\
0&0& 0
\end{bmatrix}\in \FF[s]^{p\times q}.$
Therefore, $P(s)=q_1(sr_1^T+r_2^T)$.
\item If $P(s)$ has a positive row  minimal index, then
$P(s)^T$ has a positive column  minimal index. 
Therefore, $P(s)^T=q_1(sr_1^T+r_2^T)$, i.e., $P(s)=(sr_1+r_2)q_1^T$, as desired.
\end{enumerate}
\hfill $\Box$

\section{Statement of the problem}
\label{secproblem}
The problem we deal with in this paper is the following:

\begin{problem}[Rank-one perturbation of  matrix pencils]\label{problem}
Given two matrix pencils  $A(s), B(s)\in \FF[s]^{p \times q}$,
 find necessary and sufficient conditions for the existence of a 
matrix pencil
 $P(s)\in \FF[s]^{p \times q}$ of $\rank( P(s))=1$ such that  $A(s)+P(s)\se B(s)$.  
\end{problem}

First of all we analyze two particular cases.
  \begin{itemize}
\item 
$p=1$ or $q=1$, and $A(s)\neq 0$ or $B(s)\neq 0$. If  $\FF\neq \{0, 1\}$,  
there always exists 
$P(s)=P_0+sP_1\in \FF[s]^{p\times q}$ of $\rank (P(s))=1$ such that $A(s)+P(s)\se P(s)$. For example, if $p=1$, let $c\in \FF \setminus \{0\}$ be such that 
$A(s)\neq cB(s)$. Then $A(s)+(cB(s)-A(s))\se B(s)$. If  $\FF= \{0, 1\}$, then there exists $P(s)\in \FF[s]^{p\times q}$ such that $\rank (P(s))=1$ and $A(s)+P(s)\se B(s)$ if and only if $A(s)\neq B(s)$.

\item
$p>1$ or $q>1$, and $A(s), B(s)\in \FF[s]^{p \times q}$  are such that 
 $A(s)\se B(s)$ and 
$A(s)\neq 0$. Then there always exists 
 $P(s)\in \FF[s]^{p \times q}$ of  $\rank( P(s))=1$ such that $A(s)+P(s)\se B(s)$.  For example, let  $q=2$ and let $a_1(s)\neq 0$, $a_2(s)$ be the columns of $A(s)$.
Then, $B(s)\se A(s)\se A(s)\begin{bmatrix}1&1\\0&1\end{bmatrix}=
A(s)+\begin{bmatrix}0& a_1(s)\end{bmatrix}$.
\end{itemize}

\medskip

The next lemma  shows that in order to solve Problem  \ref{problem} the pencil  $A(s)$ can  be substituted by any other pencil strictly equivalent to $A(s)$. It was proven in \cite[Lemma 3.2]{BaRo18} for
$p=q$. The proof for the general case is completely analogous.

\begin{lemma}\label{lemmasust}
 Let  $A(s), B(s), P(s)\in \FF[s]^{p\times q}$ be matrix pencils.
Let $Q\in \Gl_p(\FF)$, $R\in \Gl_q(\FF)$ and $ A'(s)=Q A(s)R$. If   $A(s)+P(s)\se B(s)$ then $A'(s)+QP(s)R\se B(s)$.
\end{lemma}
   
Problem \ref{problem} can be stated as a  pencil completion problem, as we see next.

\begin{lemma}
\label{lemmaeq}
Let  $A(s), B(s)\in \FF[s]^{p \times q}$ be matrix pencils such that 
$A(s)\not \se B(s)$.
Then there exists a matrix pencil  $P(s)\in \FF[s]^{p\times q}$ of  $\rank P(s)=1$ such that
$A(s)+P(s)\se B(s)$  if and only if one of the following conditions holds:
\begin{enumerate}
\item[(i)]
There exist matrix pencils  $a(s), b(s)\in\FF[s]^{1\times q}$, $A_{21}(s)\in \FF[s]^{(p-1)\times q}$ such that
$ A(s)\se\begin{bmatrix}a(s)\\A_{21}(s)\end{bmatrix}$ and
$B(s)\se\begin{bmatrix}b(s)\\A_{21}(s)\end{bmatrix}$.
\item[(ii)]
There exist matrix pencils $\bar a(s), \bar b(s)\in\FF[s]^{p\times 1}$, $A_{12}(s)\in \FF[s]^{p\times (q-1)}$ such that
$ A(s)\se\begin{bmatrix}\bar a(s)&A_{12}(s)\end{bmatrix}$ and
$B(s)\se\begin{bmatrix}\bar b(s)&A_{12}(s)\end{bmatrix}$.
\end{enumerate}
\end{lemma}

{\it Proof.}
Assume that  there exists a matrix pencil  $P(s)\in \FF[s]^{p\times q}$ of  $\rank P(s)=1$ such that
$A(s)+P(s)\se B(s)$.
By  Proposition \ref{pencil1}, there exist nonzero pencils $ u(s)\in \FF[s]^{p}$, $v(s)\in \FF[s]^{q}$ such that 
$P(s)=u(s)v(s)^T$ and $u(s)=u\in \FF^{p}$ or  $v(s)= v\in \FF^{q}$.
 
If  $u(s)=u\in \FF^{p}$, let  $R\in \Gl(p)$ be such that  $Ru=\begin{bmatrix}1\\0\end{bmatrix}\in\FF^{(1+(p-1))}$ and let
 $RA(s)=\begin{bmatrix}a(s)\\A_{21}(s)\end{bmatrix}\in \FF[s]^{(1+(p-1))\times q}$ and $b(s)=a(s)+v(s)^T$.
Then $A(s)\se \begin{bmatrix}a(s)\\A_{21}(s)\end{bmatrix}$ and
$B(s)\se R(A(s)+P(s))=\begin{bmatrix}a(s)\\A_{21}(s)\end{bmatrix}+\begin{bmatrix}v(s)^T\\0\end{bmatrix}=\begin{bmatrix}b(s)\\A_{21}(s)\end{bmatrix}$. 
Therefore, (i) holds.
 
If $v(s)= v\in \FF^{q}$, we can analogously obtain (ii).

\medskip

Conversely, let us assume that (i) holds.
As $A(s)\not \se B(s)$, we have $a(s)\neq b(s)$. Let  $\bar P(s)=\begin{bmatrix}b(s)-a(s)\\0\end{bmatrix}\in \FF[s]^{(1+(p-1))\times q}$.
Then $\rank \bar P(s)=1$ and $\begin{bmatrix}a(s)\\A_{21}(s)\end{bmatrix}+\bar P(s)= \begin{bmatrix}b(s)\\A_{21}(s)\end{bmatrix}$.
By  Lemma \ref{lemmasust}, there exists a pencil $P(s)$ such that  $\rank P(s)=1$ and   $A(s)+P(s)\se B(s)$.

If (ii) holds,  the result holds applying the previous case to $A(s)^T$ and
$B(s)^T$.
  
\hfill $\Box$

\section{Matrix pencil completion theorems} \label{seccompletion} 

According to Lemma \ref{lemmaeq}, the Problem \ref{problem} 
can be approached as a matrix pencil completion problem. 
We introduce in this section some results in that area which will be used later.
To state them we need some notation and  definitions.

\medskip

Given two integers $n$ and $m$, whenever $n>m$ we take  $\sum_{i=n}^{m}=0$. In the same way, if a condition is stated for $n\leq i\leq m$ with $n>m$, we understand that the condition disappears.

\medskip

 Given a  sequence of integers $a_1, \dots, a_m$ such that $a_1\geq \dots \geq a_m$ we will write $\ba=(a_1, \dots,  a_m)$ and  we will take $a_i=\infty$ for $i<1$ and $a_i=-\infty$ for $i>m$. If $a_m\geq 0$, the sequence $\ba=(a_1, \dots,  a_m)$ is called a {\em partition}.

Given three sequences $\bd$, $\ba$ and $\bg$, we introduce next the concept of generalized majorization.
\begin{definition}[\mbox{Generalized majorization \cite[Definition 2]{DoStEJC10}}]
Let $\bd = (d_1, \dots, d_m)$, $\ba=(a_1, \dots, a_s)$ and $\bg=(g_1, \dots, g_{m+s})$  be sequences of  integers
such that $d_1 \geq \dots \geq d_m$, $a_1 \geq \dots \geq a_s$, $g_1 \geq \dots \geq g_{m+s}$. We say that  $\bg$ is majorized by $\bd$ and $\ba$  $(\bg \prec' (\bd,\ba))$ if
\begin{equation}\label{gmaj1}
d_i\geq g_{i+s}, \quad 1\leq i\leq m,
\end{equation}
\begin{equation}\label{gmaj2}
\sum_{i=1}^{h_j}g_i-\sum_{i=1}^{h_j-j}d_i\leq \sum_{i=1}^j a_i, \quad 1\leq j\leq s,
\end{equation}
where $h_j=\min\{i\; : \; d_{i-j+1}<g_i\}, \   1\leq j\leq s$ \  $(d_{m+1}=-\infty)$,
\begin{equation}\label{gmaj3}
\sum_{i=1}^{m+s}g_i=\sum_{i=1}^md_i+\sum_{i=1}^sa_i.
\end{equation}
\end{definition}

In the case that $s=0$, condition (\ref{gmaj2}) disappears, and conditions (\ref{gmaj1}) and (\ref{gmaj3}) are equivalent to  $\bd=\bg$. 

\medskip

In the case that $s=1$, from condition (\ref{gmaj3}) we observe  that 
$a_1$ is completely determined by $\bd$ and $\bg$ ($a_1=\sum_{i=1}^{m+1}g_i-\sum_{i=1}^{m}d_i$), therefore we will write 
 $\bg\prec'(\bd, a_1)$ as $\bg\prec'\bd$ and we will refer to it as {\em $1$step-generalized majorization}. Moreover, it is easy to see that $\bg\prec'\bd$ if and only if
$$d_{i}=g_{i+1}, \quad h \leq i \leq m, $$
where $h=\min\{i: d_i<g_i\}$.

\begin{remark} \label{aux} \
\begin{enumerate}
\item \label{aux1} If $\bg$ and $\bd$ satisfy that  $g_i\leq d_i$ for $1\leq i \leq m$, then  $h=m+1$ and  $\bg\prec'\bd$. 
\item \label{aux2} Notice that if  $\bg\prec'\bd$ and for some index $1\leq i \leq m$ we have   $d_i>g_{i+1}$, then $i<h$. 
\item \label{aux3} In \cite{DoStEJC10}, the 1step-generalized majorization is called {\em elementary generalized majorization} and it is denoted by $\bg \prec'_1 (\bd, a_1)$.
\end{enumerate}
\end{remark}

\medskip

Given two pencils 
$H_1(s)\in \FF[s]^{(n+p) \times (n+m)}$ and  $H(s)\in \FF[s]^{ (n+p+x+y)\times (n+m)}$, of  $\rank(H_1(s))=n$ and $\rank(H(s))=n+x$, in \cite[Theorem 4.3]{DoSt19} (see also \cite[Theorem 2]{Do13}), necessary and sufficient conditions 
are given  for the existence of a pencil $Y(s)\in \FF[s]^{ (x+y)\times (n+m)}$ such that 
$H(s)\se \begin{bmatrix}
  H_1(s)\\Y(s)\end{bmatrix}$. 
The two following lemmas are particular cases of \cite[Theorem 4.3]{DoSt19} for
 $x+y=1$. 
First, we state  the result when  $x=0$, $y=1$.

\begin{lemma}[\mbox{ \cite[Particular case of  Theorem 4.3]{DoSt19}}]
      \label{lemmaDox0}
Let  $H_1(s)\in \FF[s]^{(n+p) \times (n+m)}$, $H(s)\in \FF[s]^{ (n+p+1)\times (n+m)}$ be matrix pencils of   $\rank(H_1(s))=\rank(H(s))=n$.

Let 
$\pi^1_1(s, t)\mid \dots \mid \pi^1_{n}(s, t)$,
$g_1 \geq \dots \geq g_{m}\geq 0$, and
$w_1 \geq \dots \geq w_{\theta}> 0=w_{\theta+1} \geq \dots \geq w_{p}$
be the homogeneous invariant factors, the column and the row minimal indices of $H_1(s)$, respectively, and let 
$\pi_1(s, t)\mid \dots \mid \pi_{n}(s, t)$, 
$c_1 \geq \dots \geq c_{m}\geq 0$, and
$u_1 \geq \dots \geq u_{\bar \theta}> 0=u_{\bar \theta+1} \geq \dots \geq u_{p+1}$
be the homogeneous invariant factors, the column and the row minimal indices  indices of $H(s)$, respectively.

Let $\bg=(g_1, \dots, g_m)$, $\bw=(w_1, \dots, w_p)$, $\bc=(c_1, \dots, c_m)$, and
$\bu=(u_1, \dots, u_{p+1})$.

There exists a pencil $h(s)\in\FF[s]^{1 \times (n+m)}$ such that $H(s)\se \begin{bmatrix}
  h(s)\\H_1(s)\end{bmatrix}$ if and only if
\begin{equation}\label{theta}
\bar \theta \geq \theta,
\end{equation}
\begin{equation}\label{inter}
\pi_i(s, t)\mid\pi^1_i(s, t)\mid\pi_{i+1}(s, t), \quad 1\leq i \leq n, 
\end{equation}
\begin{equation}\label{row}
\bu\prec' \bw,
\end{equation}
\begin{equation}\label{colequal}
\bg=\bc.
\end{equation}
\end{lemma}

Next, we state  the result when  $x=1$, $y=0$.

\begin{lemma}
[\mbox{ \cite[Particular case of  Theorem 4.3]{DoSt19}}]
      \label{lemmaDox1}
Let  $H_1(s)\in \FF[s]^{(n+p) \times (n+m)}$, $H(s)\in \FF[s]^{ (n+p+1)\times (n+m)}$ be matrix pencils of  $\rank(H_1(s))=n$, $\rank(H(s))=n+1$.

Let 
$\pi^1_1(s, t)\mid \dots \mid \pi^1_{n}(s, t)$,
$g_1 \geq \dots \geq g_{m}\geq 0$, and
$w_1 \geq \dots   \geq w_{p}$
be the homogeneous invariant factors, the column and the   row minimal indices of $H_1(s)$,  respectively, and let 
$\pi_1(s, t)\mid \dots \mid \pi_{n+1}(s, t)$,  
$c_1 \geq \dots \geq c_{m-1}\geq 0$, and 
$u_1 \geq \dots  \geq u_{p}$
be the homogeneous invariant factors, the column and the   row minimal indices of 
$H(s)$, respectively.

Let $\bg=(g_1, \dots, g_m)$, $\bw=(w_1, \dots, w_p)$, $\bc=(c_1, \dots, c_{m-1})$, and
$\bu=(u_1, \dots, u_{p})$.

There exists a pencil 
 $h(s)\in\FF[s]^{1 \times (n+m)}$ such that $H(s)\se \begin{bmatrix}
  h(s)\\H_1(s)\end{bmatrix}$ if and only if
(\ref{inter}), 
\begin{equation}\label{col1}
\bg\prec' \bc,
\end{equation}
\begin{equation}\label{rowequal1}
\bw=\bu.
\end{equation}

\end{lemma}

\section{Technical results}
\label{sectechnical}

The next results are  technical lemmas about 1step-generalized majorizations.
All of the sequences they deal with  are ordered partitions of nonnegative integers.

\begin{lemma}\label{propbat}
Let  $S\geq 0$ be a nonnegative integer and let  $\ba=(a_1, \dots,  a_m)$ be a partition  of nonnegative integers.  Then there exists a partition  of nonnegative integers 
$\bg=(g_1, \dots,  g_{m+1})$  such that $\sum_{i=1}^{m+1} g_i=S$ and $\bg\prec'\ba$.
\end{lemma}

{\it Proof.}
Put $a_0=+\infty$, $a_{m+1}=-\infty$. Then,
$$
\sum_{j=i}^ma_j+ia_{i-1}\geq \sum_{j=i+1}^ma_j+(i+1)a_{i}, \quad 1\leq i \leq m+1,
$$
and $
S\geq (m+2)a_{m+1}
$.
Let 
$
k=\min\{i \in \{1, \dots, m+1\}\; : \; S\geq \sum_{j=i+1}^ma_j+(i+1)a_{i}\}
$,
i.e.,
$$\sum_{j=k}^ma_j+ka_{k-1}>S\geq \sum_{j=k+1}^ma_j+(k+1)a_{k}=\sum_{j=k}^ma_j+ka_{k}.$$
Let $S'=S-\sum_{j=k}^ma_j$. Then 
$ka_{k-1}> S'\geq ka_{k}$.
Let $q$ and $r$ be the quotient and the remainder of the euclidian division of $S'$ by $k$, i.e.,  $S'=kq+r$ with $0\leq r <k$. Then
$a_{k-1}> q\geq a_{k}$.
Observe that if $k\leq m$, then $a_k\geq 0$, and if $k=m+1$, then $S'=S\geq 0$, hence $q\geq 0$.

Let us define
$$
\begin{array}{ll}
g_i=q+1, & 1\leq i \leq r,\\
g_i=q, & r+1\leq i \leq k,\\
g_i=a_{i-1}, & k+1\leq i \leq m+1. \\
\end{array}
$$
Then $g_1\geq\dots\geq g_{m+1}\geq 0$, 
$\sum_{i=1}^{m+1} g_i=S$ and
\begin{equation}\label{ga}
\begin{array}{ll}
g_{i}\leq q+1\leq a_{k-1}\leq a_i, & 1\leq i \leq k-1,\\
g_{i+1}=a_{i}, & k\leq i \leq m.\\
\end{array}
\end{equation}
 Let  $h=\min\{i\; : \; a_i<g_i\}$. From (\ref{ga}) we derive that $h\geq k$ and $a_{i}=g_{i+1}$,  $h\leq i \leq m$.
Therefore, $(g_1,\dots,g_{m+1})\prec'\ba$.

\hfill $\Box$

\begin{example}
Given the partition $\ba=(8, 6, 5, 5, 5, 3, 1)$ ($m=7$), we show some examples of the previous result for different values of $S$.
\begin{enumerate}
\item $S=50$. Then
$\sum_{j=1}^7a_j+a_{0}>S\geq
\sum_{j=2}^7a_j+2a_{1}$. $k=1$, $S'=S-\sum_{j=1}^7a_j=17$, $q=17$, $r=0$.
It results $\bg=(17, 8, 6, 5, 5, 5, 3, 1)$.
\item $S=34$. Then
$\sum_{j=3}^7a_j+3a_{2}>S\geq
\sum_{j=4}^7a_j+4a_{3}$. $k=3$, $S'=S-\sum_{j=3}^7a_j=15$, $q=5$, $r=0$.
It results $\bg=(5, 5, 5, 5, 5, 5, 3, 1)$.
\item $S=5$. Then
$\sum_{j=8}^7a_j+8a_{7}>S\geq
\sum_{j=9}^7a_j+9a_{8}$. $k=8$, $S'=S-\sum_{j=8}^7a_j=S=5$, $q=0$, $r=5$.
It results $\bg=(1, 1, 1, 1, 1, 0, 0, 0)$.
\end{enumerate}
\end{example}

\begin{lemma}\label{propbat2}
  Let  $E\geq 0$ be a nonnegative integer and let  $\ba=(a_1, \dots,  a_m)$ be a partition of nonnegative integers.  Then there exists a partition of nonnegative integers
$\be=(e_1, \dots,  e_{m-1})$  such that $\sum_{i=1}^{m-1} e_i=E$ and $\ba\prec' \be$
 if and only if
 \begin{equation}\label{eqbat2}
 E=\sum_{i=2}^{m}a_i \mbox{ or } E\geq a_1+\sum_{i=3}^{m}a_i.
  \end{equation}
\end{lemma}

{\it Proof.}
Assume that there exists 
$\be=(e_1, \dots,  e_{m-1})$  such that $\sum_{i=1}^{m-1} e_i=E$ and $\ba\prec' \be$. 
Then
$e_i\geq a_{i+1},$ $1\leq i \leq m-1$, hence  $E\geq \sum_{i=2}^{m}a_i.$

Let us suppose that $E> \sum_{i=2}^{m}a_i.$ Then $\sum_{i=1}^{m-1}(e_i-a_{i+1})>0$, hence  there exists  $k \in\{1, \dots , m-1\}$ such that  
$e_k>a_{k+1}$. Thus $k<h=\min\{i \; : e_i<a_i\}$, which means $e_1\geq a_1$ and therefore
$$
 E=e_1+\sum_{i=2}^{m-1}e_i\geq a_1+\sum_{i=3}^{m}a_i.
$$

Conversely, let us assume that (\ref{eqbat2}) holds.
We define $e_i=a_{i+1}$,  $2\leq i \leq m-1$ and 
$e_1=E-\sum_{i=3}^{m}a_i.$ Then $\sum_{i=1}^{m-1} e_i=E$.
If  $E=\sum_{i= 2}^{m} a_i$, then $e_1=a_2$ and $\ba\prec'(e_1, \dots, e_{m-1})$.
If  $E\geq a_1+\sum_{i= 3}^{m} a_i$, then  $e_1\geq a_1\geq a_3=e_2$. 
Thus, $e_1 \geq  \dots \geq e_{m-1} \geq 0$,  $h=\min\{i \; : e_i<a_i\}>1$,
and we have  $e_i=a_{i+1}$,  $h\leq i \leq m-1$. 
Therefore, $\ba\prec' (e_1, \dots, e_{m-1})$.

\hfill $\Box$ 

Given two pairs of nonincreasing sequences of integers, $(\bd, \ba)$ and
$(\bc, \bb)$, in \cite[Theorem 5.1]{DoSt13}  the authors  solved the problem of obtaining necessary and sufficient conditions
for the existence of a sequence $\bg$ that is majorized 
(in the sense of generalized majorization) by both pairs.
The conditions are very involved. 
In the first item of the next Lemma we solve the same problem for the 1step-generalized majorization of partitions. The characterization obtained is  much more simple in this case.

\begin{lemma}\label{lemmapart}
Let  $S, E\geq 0$ be nonnegative integers and let $\bc=(c_1, \dots,  c_m)$, $\bd=(d_1, \dots,  d_m)$ be partitions of nonnegative integers such that  $\bc\neq \bd$.

Let 
 $\ell=\max\{i \; :\; c_i\neq d_i\}$,
$f=\max\{i\in\{1, \dots, \ell\}\; : \; c_i<d_{i-1}\}$ $ (d_0=+\infty),$
and
$f'=\max\{i\in\{1, \dots, \ell\}\; : \; d_i<c_{i-1}\}$  $ (c_0=+\infty).$
\begin{enumerate}
\item \label{itS}
There exists a partition 
$\bg=(g_1, \dots,  g_{m+1})$  of nonnegative integers such that $\sum_{i=1}^{m+1} g_i=S$, $\bg\prec'\bc$ and $\bg\prec'\bd$ if and only if
\begin{equation}\label{eqGcd}
  S\leq \sum_{i=1}^{m}\min\{c_i, d_i\}+ \max\{c_f, d_{f'}\}.
  \end{equation}
\item \label{itEg1}
If  $f>1$ and $f'>1$, there exists a partition of nonnegative integers
$\be=(e_1, \dots,  e_{m-1})$ such that   $\sum_{i=1}^{m-1} e_i=E$, $\bc\prec' \be$ and $\bd\prec' \be$ if and only if
\begin{equation}\label{eff'm1}
E\geq \sum_{i=1}^{m}\max\{c_i, d_i\}- \max\{c_f, d_{f'}\}.
\end{equation}
\item \label{itE1}
If   $f=1$ or $f'=1$,
 there exists a partition of nonnegative integers
$\be=(e_1, \dots,  e_{m-1})$  such that $\sum_{i=1}^{m-1} e_i=E$, $\bc\prec' \be$ and $\bd\prec' \be$ if and only if
\begin{equation}\label{eff'=1}
\begin{array}{l}
 E= \sum_{i=1}^{m}\max\{c_i, d_i\}- \max\{c_f, d_{f'}\},\\\mbox{or}\\
E\geq \sum_{i=1}^{m}\max\{c_i, d_i\}- \max\{c_{f+1}, d_{f'+1}\}.\end{array}
\end{equation}
Equivalently,
 $$E= \sum_{i=2}^{m}\max\{c_i, d_i\}\mbox{ or }
E\geq \max\{c_1, d_1\}+\sum_{i=3}^{m}\max\{c_i, d_i\}.
$$
\end{enumerate}
\end{lemma}
  
{\it Proof.}
Let us assume that $c_\ell>d_\ell$. If $d_\ell>c_\ell$ the proof is analogous.

We have $c_{\ell-1}\geq c_\ell>d_\ell,$ hence $f'=\ell$. 
Moreover, 
$c_i\geq c_{i+1}\geq d_i,$ $f\leq i \leq \ell-1.$
Then, $c_f \geq d_f\geq d_\ell=d_{f'}$.
Hence, when   $c_\ell>d_\ell$, conditions (\ref{eqGcd}), (\ref{eff'm1}) and (\ref{eff'=1}) are respectively equivalent to
\begin{equation}\label{eqGcdbis}
  S\leq \sum_{i=1}^{f-1}\min\{c_i, d_i\}+ c_f+ \sum_{i=f}^{m}d_i,
  \end{equation}
\begin{equation}\label{eff'm1bis}
E\geq \sum_{i=1}^{f-1}\max\{c_i, d_i\}+ \sum_{i=f+1}^{m}c_i,
\end{equation}
\begin{equation}\label{eff'=1bis}
 E= \sum_{i=2}^{m}c_i\mbox{ or }
E\geq c_1+\sum_{i=3}^{m}c_i.
\end{equation}
  
Moreover, if $f'=1$ then $\ell=1$ and as a consequence, $f=1$. 
Therefore, when   $c_\ell>d_\ell$, the condition $f=1$ or $f'=1$ is equivalent to $f=1$.

Let us prove the different cases.

\begin{enumerate}
\item
Assume that 
there exists a partition 
$\bg=(g_1, \dots,  g_{m+1})$  such that $\sum_{i=1}^{m+1} g_i=S$, $\bg\prec'\bc$ and $\bg\prec'\bd$.

Let 
$h=\min\{i \; : c_i<g_i\}$ and $h'=\min\{i \; : d_i<g_i\}$.
As $g_{\ell+1}\leq d_\ell<c_\ell$, by Remark \ref{aux}, item \ref{aux2}, we have  $\ell < h$. In the same way, as  $g_{f}\leq c_f<d_{f-1}$, $f-1<h'$.

Therefore,
$$
\begin{array}{ll}
g_i\leq\min\{c_i, d_i\}, & 1\leq i \leq f-1,\\
g_f\leq c_f, \\
g_{i+1}\leq d_{i}, &f\leq i \leq m,
\end{array}
$$
from where we obtain (\ref{eqGcdbis}).

Conversely, let us assume that (\ref{eqGcdbis}) holds.

\begin{itemize}
\item 
If $S<\sum_{i=1}^{f-1}\min\{c_i, d_i\}$ then $f>1$. Let
$$k=\min\{i\in\{1, \dots f-1\}\; : \; S<\sum_{i=1}^{k}\min\{c_i, d_i\}\},$$ i.e., 
$\sum_{i=1}^{k-1}\min\{c_i, d_i\}\leq S<\sum_{i=1}^{k}\min\{c_i, d_i\}$ and define
$$
\begin{array}{ll}
g_i=\min\{c_i, d_i\}, & 1\leq i \leq k-1,\\
g_k=S-\sum_{i=1}^{k-1}\min\{c_i, d_i\}, \\
g_{i}=0, &k+1\leq i \leq m+1.
\end{array}
$$
Then $\sum_{i=1}^{m+1} g_i=S$ and   $g_k<\min\{c_k, d_k\}$. Therefore, $g_1\geq \dots \geq g_{k-1}>g_k\geq 0=g_{k+1}= \dots =g_{m+1}$.
Thus, $\bg=(g_1, \dots, g_{m+1})$ is a  partition. As $g_i \leq \min\{c_i, d_i\}$, $1\leq i \leq m$, by Remark \ref{aux}, item \ref{aux1}, we have  $\bg\prec'\bc$ and $\bg\prec'\bd$.

\item
If $S\geq\sum_{i=1}^{f-1}\min\{c_i, d_i\}$, let $\overline S=S- \sum_{i=1}^{f-1}\min\{c_i, d_i\}\geq 0$. 
Then  $\overline S\leq c_f+\sum_{i=f}^{m}d_i.$
We define
$
\overline d_i=d_{f-1+i}$, $1\leq i \leq m-f+1,$ and $\overline \bd=(\overline d_1, \dots, \overline d_{m-f+1})$, i.e., 
$\overline \bd=(d_f, \dots, d_{m})$. 
By Lemma 
\ref{propbat}, there exists a partition  $\overline \bg=(\overline g_1, \dots, \overline g_{m-f+2})$ such that 
$\sum_{i=1}^{m-f+2}\overline g_i=\overline S$ and  $\overline \bg\prec' \overline \bd$.

Now we define
$$
\begin{array}{ll}
g_i=\min\{c_i, d_i\}, & 1\leq i \leq f-1,\\
g_i= \overline g_{i-f+1}, & f\leq i \leq m+1.
\end{array}
$$
Let us see that $g_{f}\leq g_{f-1}$, i.e., that $\overline g_1 \leq \min\{c_{f-1}, d_{f-1}\}$.
If  $\overline g_1 \leq \overline d_1$, then $\overline g_1 \leq d_f=\min\{c_{f}, d_{f}\}\leq \min\{c_{f-1}, d_{f-1}\}$.
If  $\overline g_1 > \overline d_1$, then $\overline d_i= \overline g_{i+1}$, $1\leq i \leq  m-f+1$, hence
$\overline S=\overline g_1+\sum_{i=1}^{m-f+1}\overline d_i$. As $\overline S\leq c_f+\sum_{i=f}^{m}d_i$, we obtain  $\overline g_1\leq c_f\leq \min\{c_{f-1}, d_{f-1}\}$.
Therefore $\bg=(g_1, \dots, g_{m+1})$ is a partition.

Let $h'=\min\{i \; : \; d_i<g_i\}$ and $\overline h'=\min\{i \; : \; \overline d_i<\overline g_i\}$. 
Observe that $d_i \geq g_i$, $1\leq i \leq f-1$, 
$d_i=\overline{d}_{i-f+1}\geq \overline{g}_{i-f+1}=g_i$
for $f\leq i < f+\overline{h}'-1$,  and 
$d_{f+\overline{h}'-1}=\overline{d}_{\overline{h}'}<\overline{g}_{\overline{h}'}= g_{f+\overline{h}'-1}$.
Then, $h'=f+\overline{h}'-1$. As $d_i=\overline d_{i-f+1}=\overline g_{i-f+2}=g_{i+1}$ for $h'\leq i \leq m$, we obtain that $\bg \prec' \bd$. 

Let $h=\min\{i \; : \; c_i<g_i\}$. Recall that $d_i\leq c_i$ for $f\leq i \leq m$. We have,
$$
\begin{array}{ll}
g_i\leq c_i, & 1\leq i \leq f-1,\\
g_i\leq d_i\leq c_i, & f\leq i \leq h'-1,\\
g_i\leq d_{i-1}\leq c_i, & f+1\leq i \leq \ell.\\
\end{array}
$$
Thus,  $h>\max\{h'-1, \ell\}$ and, as a consequence,  
 $c_i=d_i=g_{i+1}$ for $h\leq i \leq m$. Therefore $\bg \prec' \bc$.
\end{itemize}

\item
Assume that $f>1$ (hence $f'>1$) and there exists a partition 
$\be=(e_1, \dots,  e_{m-1})$ such that   $\sum_{i=1}^{m-1} e_i=E$, $\bc\prec' \be$ and $\bd\prec' \be$.
Then 
$e_i\geq c_{i+1},$ $1\leq i \leq m-1.$
Moreover, 
$e_{\ell-1}\geq c_\ell>d_\ell$, by Remark \ref{aux}, item \ref{aux2}, 
$e_{i}\geq d_i,$  $1\leq i \leq \ell-1$.
Hence,
$e_{f-1}\geq d_{f-1}>c_f$, and as before it means that $e_{i}\geq c_{i},$ $1\leq i \leq f-1.$ 
Thus, 
$$
\begin{array}{ll}
e_i\geq \max\{c_i,d_i\} & 1\leq i \leq f-1,\\
e_i\geq  c_{i+1}, & f\leq i \leq m-1,\\
\end{array}
$$
and we obtain (\ref{eff'm1bis}).

Conversely, let us assume that $f>1$ and (\ref{eff'm1bis}) holds.
Let us define
$$
\begin{array}{ll}
e_i=\max\{c_i, d_i\}, & 2\leq i\leq f-1,\\
e_i=c_{i+1}, & f\leq i\leq m-1,\\
e_1=E-\sum_{i=2}^{m-1}e_i.
\end{array}
$$
Then,  $E=\sum_{i=1}^{m-1}e_i$,
$e_2\geq \dots \geq e_{f-1}\geq c_{f-1} \geq c_{f+1} =e_f\geq \dots \geq e_{m-1}$ and,
from (\ref{eff'm1bis}) we derive
$ e_1\geq \max\{c_1, d_1\}\geq e_2$. Therefore
$(e_1, \dots,  e_{m-1})$ is a partition.
Let $h=\min\{i \; : e_i<c_i\}$ and  $h'=\min\{i \; : e_i<d_i\}$.
It is clear that $h\geq f$ and $h'\geq f$, hence  $e_i=c_{i+1},$  $h\leq i\leq m-1$, which  means that 
$\bc\prec' \be$. 
Moreover, for $f \leq i \leq \ell-1$, $e_i=c_{i+1}\geq d_i$, thus $h'\geq \ell$
and $e_i=c_{i+1}=d_{i+1}$, $h'\leq i\leq m-1.$ Therefore $\bd\prec' \be$. 

\item
Let us assume that $f=1$
 and there exists a partition 
$\be=(e_1, \dots,  e_{m-1})$ such that   $\sum_{i=1}^{m-1} e_i=E$, $\bc\prec' \be$ and $\bd\prec' \be$.
From Lemma \ref{propbat2} we obtain (\ref{eff'=1bis}).

Conversely, let us assume that $f=1$ and
 (\ref{eff'=1bis}) holds.
By  Lemma \ref{propbat2}, 
there exists a partition 
 $\be=(e_1, \dots,  e_{m-1})$  such that  $\sum_{i=1}^{m-1} e_i=E$ and $\bc\prec' \be$.
Therefore, $e_i\geq c_{i+1}, \ 1\leq i\leq m-1$.

 Let  $h=\min\{i \; : e_i<c_i\}$ and $h'=\min\{i \; : e_i<d_i\}$. 
 As $f=1$, we have $ c_{i}\geq d_{i-1}\geq d_i$,  $2\leq i \leq  \ell$. Therefore, $e_i\geq c_{i+1}\geq d_i$, $1\leq i<\ell$ and $h'\geq \ell$.
Since $e_{h'}<d_{h'}=c_{h'}$,  $h\leq h'$ 
and $e_{i}=c_{i+1}=d_{i+1}$, $h'\leq i \leq m-1$.
Hence,  $\bd\prec' \be$.

\end{enumerate}

\hfill $\Box$

\begin{remark}\label{remarkL54}
Observe that condition (\ref{eff'=1}) implies condition (\ref{eff'm1}).
\end{remark}

\begin{lemma}
\label{lemmacard}
Let $\ba=(a_1, \dots,  a_{m})$, $\be=(e_1, \dots,  e_{m-1})$ be partitions of nonnegative integers such that   $\ba\prec' \be$ and 
$
 \sum_{i=1}^{m-1}e_i \leq \sum_{i=1}^{m}a_i.
$
Let $\theta=\#\{i \; : \; e_i>0\}$ and  $\overline \theta=\#\{i \; : \; a_i>0\}$. Then $\overline \theta \geq \theta$.
\end{lemma}

{\it Proof.}
We have  
$\theta\leq m-1$,  $\overline \theta \leq m$ and  $
 \sum_{i=1}^{\theta}e_i \leq \sum_{i=1}^{\overline \theta}a_i.$
Let $h=\min \{i \; : \; e_i<a_i\}$. Then $e_{i}=a_{i+1}$ for  $h\leq i \leq m-1$.

Assume that  $\theta>\overline \theta$. Then
$0<\sum_{i=\overline \theta+1}^{\theta}e_i \leq \sum_{i=1}^{\overline \theta}(a_i-e_i).$
It means that  there exists  $i \in \{1,\dots, \overline \theta  \}$ such that $a_i-e_i>0$. Therefore, $h\leq \overline \theta <\theta <m$, from where we conclude that 
$e_{\overline \theta}=a_{\overline \theta+1}=0$, which is a contradiction with   $\theta>\overline \theta$.

\hfill $\Box$

\section{Main theorem}
\label{secmain}

In the following Theorem we give a solution to Problem \ref{problem}.

\begin{theorem}\label{maintheogen}
Let $A(s), B(s)\in \FF[s]^{p\times q}$ be matrix pencils such that  $A(s)\not \se B(s)$. 
Let $\rank A(s)=n_1$, $\rank B(s)=n_2$, let 
$\phi_1(s, t)\mid \dots \mid \phi_{n_1}(s, t)$,
$c_1 \geq \dots \geq c_{q-n_1}\geq  0$ and
$u_1 \geq \dots \geq u_{p-n_1}\geq  0$
 be, respectively,  the homogeneous invariant factors, column minimal indices and row minimal indices of $A(s)$ and let 
$\psi_1(s, t)\mid \dots \mid \psi_{n_2}(s, t)$, 
$d_1 \geq \dots \geq d_{q-n_2}\geq  0$ 
and 
$v_1 \geq \dots \geq v_{p-n_2}\geq  0$
 be, respectively,  the homogeneous invariant factors, column minimal indices and row minimal indices of $B(s)$.

Let
$n=\min\{n_1, n_2\}$, 
$\bc=(c_1, \dots, c_{q-n_1})$, $\bd=(d_1, \dots, d_{q-n_2})$,
$\bu=(u_1, \dots, u_{p-n_1})$ and $\bv=(v_1, \dots, v_{p-n_2})$.

\begin{enumerate}
\item \label{equal}
If $\bc=\bd$, $\bu=\bv$, then there exists a pencil $P(s)\in \FF[s]^{p \times q}$ of    $\rank (P(s))=1$ such that 
$A(s)+P(s)\se B(s)$ if and only if
\begin{equation}\label{eqintfihr1}
\psi_{i-1}(s, t)\mid\phi_i(s, t)\mid\psi_{i+1}(s, t), \quad 1\leq i \leq n.
\end{equation}
\item \label{rowequal}
If $\bc\neq \bd$, $\bu=\bv$,
 let
$$\ell=\max\{i \; :\; c_i\neq d_i\},$$
$$
f=\max\{i\in\{1, \dots, \ell\}\; : \; c_i<d_{i-1}\} \quad (d_0=+\infty),
$$
$$
f'=\max\{i\in\{1, \dots, \ell\}\; : \; d_i<c_{i-1}\} \quad (c_0=+\infty),
$$
$$
G=n-1-\sum_{i=1}^{n-1}\deg(\gcd(\phi_{i+1}(s, t),\psi_{i+1}(s, t)))- \sum_{i=1}^{p-n}u_i,
$$
$$
T=n-\sum_{i=1}^{n}\deg(\lcm(\phi_{i}(s, t),\psi_{i}(s, t)))- \sum_{i=1}^{p-n}u_i.
$$
Then there exists a pencil $P(s)\in \FF[s]^{p \times q}$ of  $\rank (P(s))=1$ such that 
$A(s)+P(s)\se B(s)$ if and only if
 (\ref{eqintfihr1}) and one of the two following conditions holds:
\begin{equation}\label{eqGg}
 G\leq \sum_{i=1}^{q-n}\min\{c_i, d_i\}+ \max\{c_f, d_{f'}\},
\end{equation}
or
\begin{equation}\label{eqT}
T\geq \sum_{i=1}^{q-n}\max\{c_i, d_i\}- \max\{c_f, d_{f'}\}.
\end{equation}
\item  
\label{colequalcdrns}
If $\bc= \bd$, $\bu\neq \bv$,
let
 $$\bar \ell=\max\{i \; :\; u_i\neq v_i\},$$
$$
\bar f=\max\{i\in\{1, \dots, \bar \ell\}\; : \; u_i<v_{i-1}\} \quad (v_0=+\infty),
$$
$$
\bar f'=\max\{i\in\{1, \dots, \bar \ell\}\; : \; v_i<u_{i-1}\} \quad (u_0=+\infty),
$$
$$
\bar G=:n-1-\sum_{i=1}^{n-1}\deg(\gcd(\phi_{i+1}(s, t),\psi_{i+1}(s, t)))-\sum_{i=1}^{q-n}c_i,
$$
$$
\bar T=n-\sum_{i=1}^{n}\deg(\lcm(\phi_{i}(s, t),\psi_{i}(s, t)))-\sum_{i=1}^{q-n}c_i.
$$
Then there exists a pencil $P(s)\in \FF[s]^{p \times q}$ of $\rank (P(s))=1$ such that 
$A(s)+P(s)\se B(s)$ if and only if
 (\ref{eqintfihr1}) and one of the two following conditions holds:
\begin{equation}\label{eqbarGg}
 \bar G\leq \sum_{i=1}^{p-n}\min\{u_i, v_i\}+ \max\{u_{\bar f}, v_{\bar f'}\},
\end{equation}
or
\begin{equation}\label{eqbarT}
\bar T\geq \sum_{i=1}^{p-n}\max\{u_i, v_i\}- \max\{u_{\bar f}, v_{\bar f'}\}.
\end{equation}
\item \label{inequal}
If $\bc\neq \bd$, $\bu\neq \bv$,
then there exists a pencil $P(s)\in \FF[s]^{p \times q}$ of $\rank (P(s))=1$ such that 
$A(s)+P(s)\se B(s)$ if and only if
there exist homogeneous polynomials 
$\pi_1^1(s,t)\mid \dots  \mid \pi_{n}^1(s,t)$ such that 
\begin{equation}\label{boca}
\lcm(\phi_i(s,t), \psi_i(s,t))\mid \pi_i^1(s,t)\mid \gcd(\phi_{i+1}(s,t), \psi_{i+1}(s,t)), \quad 1\leq i \leq n.
\end{equation}
and  one of the four following conditions holds:
\begin{enumerate}
\item[(a)]
\begin{equation}\label{cdsr}
\bc\prec'\bd, \quad \bu\prec'\bv,
\end{equation}
\begin{equation}\label{sumpi1}
\sum_{i=1}^{n} \deg(\pi_i^1(s,t))=n-\sum_{i=1}^{q-n_1}c_i-\sum_{i=1}^{p-n_2}v_i.
\end{equation}
\item[(b)]
\begin{equation}\label{dcrs}
\bd\prec'\bc, \quad \bv\prec'\bu,
\end{equation}
\begin{equation}\label{sumpi1t}
\sum_{i=1}^{n} \deg(\pi_i^1(s,t))=n-\sum_{i=1}^{q-n_2}d_i-\sum_{i=1}^{p-n_1}u_i.
\end{equation}

\item[(c)]
  (\ref{cdsr}) and (\ref{sumpi1t}).
 \item[(d)]
  (\ref{dcrs}) and (\ref{sumpi1}). 

\end{enumerate}

\end{enumerate}

\end{theorem}

{\it Proof.}
\underline{Necessity}. Let us assume that there exists a pencil
  $P(s)\in \FF[s]^{p\times q}$ of   $\rank P(s)=1$ such that 
$A(s)+P(s)\se B(s)$. By Lemma \ref{lemmaeq},  one of the two following conditions holds:
\begin{enumerate}
\item[(i)]
There exist pencils $a(s), b(s)\in\FF[s]^{1\times q}$ and $A_{21}(s)\in \FF^{(p-1)\times q}$ such that 
$ A(s)\se\begin{bmatrix}a(s)\\A_{21}(s)\end{bmatrix}$ and
$B(s)\se\begin{bmatrix}b(s)\\A_{21}(s)\end{bmatrix}$.
\item[(ii)]
There exist pencils $\bar a(s), \bar b(s)\in\FF[s]^{p\times 1}$ and $A_{12}(s)\in \FF^{p\times (q-1)}$ such that 
$ A(s)\se\begin{bmatrix}\bar a(s)&A_{12}(s)\end{bmatrix}$ and 
$B(s)\se\begin{bmatrix}\bar b(s)&A_{12}(s)\end{bmatrix}$. 
\end{enumerate}

\begin{itemize}
\item
Let us assume that (i) holds.
Then  $n\geq \rank(A_{21}(s))\geq \max\{n_1, n_2\}-1\geq n-1$, hence
 $\rank(A_{21}(s))=n-x$ with $x=0$ or $x=1$.
Let 
$\pi_i^1(s, t)\mid \dots \mid \pi^1_{n-x}(s,t)$, $\bg=(g_1, \dots, g_{q-n+x})$ and $\bw=(w_1, \dots  w_{p-1-n+x})$ be, respectively,  the homogeneous invariant factors , column minimal indices and row minimal indices of
$A_{21}(s)$.
By Lemmas \ref{lemmaDox0} and \ref{lemmaDox1},
\begin{equation}\label{interx}
\begin{array}{ll}
  \phi_i(s, t)\mid\pi^1_i(s, t)\mid\phi_{i+1}(s, t),&
1\leq i \leq n-x,
\\
\psi_i(s, t)\mid\pi^1_i(s, t)\mid\psi_{i+1}(s, t),
  & 1\leq i \leq n-x.
\end{array}
\end{equation}
Thus,
$$
\begin{array}{ll}
  \psi_{i-1}(s, t)\mid\phi_{i}(s, t),&
1\leq i \leq n,
\\
\phi_{i}(s, t)\mid\psi_{i+1}(s, t),
  & 1\leq i \leq n-x.
\end{array}
$$

Notice that in the case that   $x=1$ we have  $n_1=n_2=n$ and $\phi_n(s,t)\mid \psi_{n+1}(s, t)=0$ is also satisfied. Therefore,  (\ref{eqintfihr1}) holds.

\begin{enumerate}
\item
 
Assume that  $\bc= \bd$, $\bu= \bv$. As  it has been seen, condition (\ref{eqintfihr1}) is necessary.

\item
Assume that  $\bc\neq \bd$, $\bu= \bv$. Then $n_1=n_2=n$. If $\rank(A_{21}(s))=n$, then  from 
Lemma \ref{lemmaDox0} we obtain $\bg=\bc$ and  $\bg=\bd$, which is a contradiction with
$\bc\neq \bd$.
Therefore,  $\rank(A_{21}(s))=n-1$, i.e., $x=1$. 
Applying Lemma \ref{lemmaDox1}, we obtain
\begin{equation}\label{colnnn-1}
\bg\prec' \bc,\quad \bg\prec' \bd,
\end{equation}
\begin{equation}\label{rowequalnn-1}
\bw=\bu=\bv.
\end{equation}
From (\ref{interx}) and (\ref{rowequalnn-1}),
$$\begin{array}{rl}\sum_{i=1}^{q-n+1}g_i=&n-1-\sum_{i=1}^{n-1}\deg(\pi_i^1(s, t)-\sum_{i=1}^{p-n}w_i\\\geq &
n-1-\sum_{i=1}^{n-1}\deg(\gcd(\phi_{i+1}(s,t), \psi_{i+1}(s,t)))-\sum_{i=1}^{p-n}u_i=G.\end{array}
$$
By Lemma \ref{lemmapart}, 
$$\sum_{i=1}^{q-n+1}g_i\leq \sum_{i=1}^{q-n}\min\{c_i, d_i\}+ \max\{c_f, d_{f'}\}.$$
Therefore, 
(\ref{eqGg}) holds.

\item
Assume that $\bc=\bd$, $\bu\neq \bv$. Then $n_1=n_2=n$. If $\rank(A_{21}(s))=n-1$,  then from 
Lemma \ref{lemmaDox1}, we obtain $\bw=\bu$ and  $\bw=\bv$, which is a contradiction with $\bu\neq \bv$.
Therefore, $\rank(A_{21}(s))=n$, i.e.,  $x=0$. 

Applying Lemma \ref{lemmaDox0}, we obtain
\begin{equation}\label{rownnn}
\bu\prec' \bw, \quad \bv\prec' \bw,
\end{equation}
\begin{equation}\label{colequalnnnT}
\bg=\bc=\bd.
\end{equation}
From (\ref{interx}) and (\ref{colequalnnnT}),
$$\begin{array}{rl}
\sum_{i=1}^{p-n-1}w_i=&n-\sum_{i=1}^{n}\deg(\pi_i^1(s, t)-\sum_{i=1}^{q-n}g_i\\\leq &
n-\sum_{i=1}^{n}\deg(\lcm(\phi_{i}(s,t), \psi_{i}(s,t)))-\sum_{i=1}^{q-n}c_i=\bar T.
\end{array}
$$

By Lemma \ref{lemmapart} and Remark \ref{remarkL54}, 
$$\sum_{i=1}^{p-n-1}w_i\geq \sum_{i=1}^{p-n}\max\{u_i, v_i\}- \max\{u_{\bar f}, v_{\bar f'}\}.$$
Therefore, (\ref{eqbarT}) holds.

\item 
Assume that $\bc\neq \bd$, $\bu\neq \bv$. If $\rank(A(s))=\rank (B(s))$, then applying 
Lemmas \ref{lemmaDox0} and \ref{lemmaDox1}, we obtain $\bg=\bc=\bd$ or $\bw=\bu=\bv$, which is a contradiction.
Therefore,  $\rank(A(s))\neq \rank (B(s))$. Then  $n\geq \rank (A_{21}(s))\geq \max\{n_1, n_2\}-1=n$, i.e., 
 $\rank (A_{21}(s))=n$ ($x=0$).
 From (\ref{interx}) we derive (\ref{boca}).

If $\rank (A(s))<\rank (B(s))$, 
then $\rank (A(s))=n$, $\rank (B(s))=n+1$.
Applying Lemmas  \ref{lemmaDox0} and \ref{lemmaDox1} we obtain
\begin{equation}\label{colnn+1n}
\bg= \bc,\quad \bu\prec'\bw,
\end{equation}
\begin{equation}\label{rownn+1n}
\bg\prec' \bd, \quad \bw=\bv.
\end{equation}
From
(\ref{colnn+1n}) and (\ref{rownn+1n}) we derive (\ref{cdsr}) and
(\ref{sumpi1}).

Analogously, if $\rank (B(s))<\rank (A(s))$ we obtain (\ref{dcrs}) and
(\ref{sumpi1t}).
  
\end{enumerate}

\item
Let us assume that  (ii) holds.
Then 
$$ A(s)^T\se\begin{bmatrix}\bar a(s)^T\\A_{12}(s)^T\end{bmatrix}, \quad 
B(s)^T\se\begin{bmatrix}\bar b(s)^T\\A_{12}(s)^T\end{bmatrix}.$$

Recall that the column and row minimal indices of a pencil are, respectively, the row and column minimal indices of its transposed.

Applying the results of the previous case, the interlacing condition (\ref{eqintfihr1}) is satisfied and

\begin{itemize}
\item
If  $\bc=\bd$, $\bu\neq \bv$  we obtain (\ref{eqbarGg}).

\item
If  $\bc\neq \bd$, $\bu= \bv$ we obtain (\ref{eqT}).

\item
If  $\bc\neq d$, $\bu\neq \bv$ and $\rank (A(s))<\rank (B(s))$ we obtain (\ref{cdsr}) and
(\ref{sumpi1t}).

\item If  $\bc\neq \bd$, $\bu\neq \bv$ and $\rank (B(s))<\rank (A(s))$
we obtain (\ref{dcrs}) and
(\ref{sumpi1}).
\end{itemize}

\end{itemize}

\underline{Sufficiency}.

\medskip

\noindent
 \underline{Case $\bu=\bv$}. In this case, $n=n_1=n_2$.
\begin{itemize}
\item
Assume that $\bc=\bd$ and (\ref{eqintfihr1}) holds or that $\bc\neq \bd$ and (\ref{eqintfihr1}) and (\ref{eqGg}) hold.
By Lemma \ref{lemmaeq}, it is enough to prove the existence of matrix pencils  $a(s), b(s)\in\FF[s]^{1\times q}$, $A_{21}(s)\in \FF[s]^{(p-1)\times q}$ such that
$ A(s)\se\begin{bmatrix}a(s)\\A_{21}(s)\end{bmatrix}$ and
$B(s)\se\begin{bmatrix}b(s)\\A_{21}(s)\end{bmatrix}$.

Let 
$$\pi^1_i(s, t)=\gcd(\phi_{i+1}(s, t),\psi_{i+1}(s, t)), \quad 1\leq i \leq n-1.$$
Then
$\pi_1^1(s,t)\mid \dots \mid \pi_{n-1}^1(s,t)$ and (\ref{eqintfihr1}) implies
that
\begin{equation}\label{intern-1}
\begin{array}{ll}
  \phi_i(s, t)\mid\pi^1_i(s, t)\mid\phi_{i+1}(s, t),&
1\leq i \leq n-1,
\\
\psi_i(s, t)\mid\pi^1_i(s, t)\mid\psi_{i+1}(s, t),
  & 1\leq i \leq n-1.
\end{array}
\end{equation}

Let $S=n-1-\sum_{i=1}^{n-1}\deg(\pi^1_i(s, t))-\sum_{i=1}^{p-n}u_i$ and
let us see that $S\geq 0$. 
 We have 
 $$\begin{array}{rl}
\sum_{i=1}^{n-1}\deg(\pi^1_i(s, t))+\sum_{i=1}^{p-n}u_i\leq & \sum_{i=2}^{n}\deg(\phi_i(s, t))+\sum_{i=1}^{p-n}u_i\\=& n-\deg(\phi_1(s,t))-\sum_{i=1}^{q-n}c_i,\\ 
\end{array}$$
and
$$\begin{array}{rl}
\sum_{i=1}^{n-1}\deg(\pi^1_i(s, t))+\sum_{i=1}^{p-n}u_i\leq& \sum_{i=2}^{n}\deg(\psi_i(s, t))+\sum_{i=1}^{p-n}u_i\\=& n-\deg(\psi_1(s,t))-\sum_{i=1}^{q-n}d_i. \end{array}$$
If $\sum_{i=1}^{n-1}\deg(\pi^1_i(s, t))+\sum_{i=1}^{p-n}u_i= n$, then $\phi_1(s, t)= \psi_1(s, t)=1$, $\sum_{i=1}^{q-n}c_i=\sum_{i=1}^{q-n}d_i=0$,
and $$\sum_{i=1}^{n-1}(\deg(\pi^1_i(s, t))-\deg(\phi_{i+1}(s, t)))= \sum_{i=1}^{n-1}(\deg(\pi^1_i(s, t))-\deg(\psi_{i+1}(s, t)))=0,$$ 
therefore, $\bc=\bd$, $\pi^1_i(s, t)=\phi_{i+1}(s, t)=\psi_{i+1}(s, t)$, $1\leq i \leq n-1$  and   $A(s)\se B(s)$. As $A(s)\not \se B(s)$,  we derive $S\geq 0$.

Notice that in the case that $\bc\neq \bd$, because of condition (\ref{eqGg}),
 $S=G\leq \sum_{i=1}^{q-n}\min\{c_i, d_i\}+ \max\{c_f, d_{f'}\}$.
Hence, by Lemma \ref{propbat} (in the case $\bc=\bd$) or by Lemma   \ref{lemmapart}  (in the case  $\bc\neq \bd$), there exists a partition of nonnegative integers
$\bg=(g_1, \dots,  g_{q-n+1})$  satisfying $\sum_{i}^{q-n+1} g_i=S$ and 
(\ref{colnnn-1}).

Since $\sum_{i=1}^{n-1}\deg(\pi^1_i(s, t))+\sum_{i=1}^{p-n}u_i+\sum_{i}^{q-n-1} g_i=n-1 $, there exists  a pencil $A_{21}(s)\in \FF^{(p-1)\times q} $ of $\rank(A_{21}(s))=n-1$,  homogeneous invariant factors $\pi^1_i(s, t)\mid \dots \mid\pi^1_{n-1}(s, t)$, column minimal indices 
$g_1\geq \dots \geq g_{q-n+1}$ and row minimal indices  $u_1\geq \dots \geq u_{p-n}$. 

From (\ref{intern-1}) and (\ref{colnnn-1}) and  Lemma \ref{lemmaDox1},   there exist pencils $a(s), b(s)\in\FF[s]^{1\times q}$,  such that
$ A(s)\se\begin{bmatrix}a(s)\\A_{21}(s)\end{bmatrix}$ and
$B(s)\se\begin{bmatrix}b(s)\\A_{21}(s)\end{bmatrix}$. 

\item 
Assume that $\bc\neq \bd$ and that (\ref{eqintfihr1}) and
 (\ref{eqT}) hold.
By Lemma \ref{lemmaeq}, it is enough to prove the existence of matrix pencils  $\bar a(s), \bar b(s)\in\FF[s]^{p\times 1}$, $A_{12}(s)\in \FF[s]^{p\times (q-1)}$ such that
$ A(s)\se\begin{bmatrix}\bar a(s)& A_{12}(s)\end{bmatrix}$ and
$B(s)\se\begin{bmatrix}\bar b(s) & A_{12}(s)\end{bmatrix}$.

Let $x=T-\sum_{i=1}^{q-n}\max\{c_i, d_i\}+\max\{c_f, d_{f'}\}$. Condition (\ref{eqT}) implies that $x\geq 0$.
Let $\gamma(s,t)$ be a homogeneous polynomial of $\deg(\gamma(s,t))= x$ and define
$$
\begin{array}{ll}
\bar \pi^1_i(s, t): =\lcm(\phi_{i}(s, t),\psi_{i}(s, t)), \quad 1\leq i \leq n-1,
\\
\bar \pi^1_n(s, t)=\gamma(s,t)\lcm(\phi_{n}(s, t),\psi_{n}(s, t)).
\end{array}
$$
Then,
$\bar \pi_1^1(s,t)\mid \dots \mid \bar \pi_{n}^1(s,t)$ and (\ref{eqintfihr1}) implies
\begin{equation}\label{intern}
\begin{array}{ll}
\phi_i(s, t)\mid\bar \pi^1_i(s, t)\mid\phi_{i+1}(s, t), & 1\leq i \leq n,\\
\psi_i(s, t)\mid\bar \pi^1_i(s, t)\mid\psi_{i+1}(s, t), & 1\leq i \leq n.
\end{array}
\end{equation}
Let $T'=T-x=\sum_{i=1}^{q-n}\max\{c_i, d_i\}- \max\{c_f, d_{f'}\}$. Then $T'\geq 0$. By  Lemma \ref{lemmapart}, there exists a partition of nonnegative integers 
$\bg=(g_1, \dots,  g_{q-n-1})$ such that $\sum_{i}^{q-n-1} g_i=T'$ and
\begin{equation}\label{revcolnnn}
\bc\prec' \bg,\quad \bd\prec' \bg.
\end{equation}

(Notice that due to the value of $T'$, the  conditions in the cases \ref{itEg1} and \ref{itE1} of Lemma \ref{lemmapart} are satisfied). 

By the definition of  $T$,
$$\begin{array}{l}
T'\leq T\leq n-\sum_{i=1}^{n}\deg(\phi_i(s,t))-\sum_{i=1}^{p-n}u_i=\sum_{i=1}^{q-n}c_i,\\
T'\leq T\leq n-\sum_{i=1}^{n}\deg(\psi_i(s,t))-\sum_{i=1}^{p-n}v_i=\sum_{i=1}^{q-n}d_i.\end{array}$$
Then,
from Lemma \ref{lemmacard} we obtain
\begin{equation}\label{eqrr}
\#\{i \; : \;\ g_i>0\}\leq \#\{i \; : \; c_i>0\}, 
\quad
\#\{i \; : \;\ g_i>0\}\leq \#\{i \; : \; d_i>0\}.
\end{equation}
As $\sum_{i=1}^{n}\deg(\bar \pi^1_i(s, t))+\sum_{i=1}^{p-n}u_i+\sum_{i}^{q-n+1} g_i=n,$ there exists  a pencil
$A_{12}(s)\in \FF^{p\times (q-1)}$ of $\rank (A_{12}(s))=n$,  homogeneous invariant factors $\bar \pi^1_i(s, t)\mid \dots \mid \bar \pi^1_{n}(s, t)$, column minimal indices 
$g_1\geq \dots \geq g_{q-n-1}$ and row minimal indices  $u_1\geq \dots \geq u_{p-n}$. 
From  (\ref{intern})-(\ref{eqrr}) and  Lemma \ref{lemmaDox0}, there exist pencils
$\bar a(s)^T, \bar b(s)^T\in\FF[s]^{1\times p}$ such that
$ A(s)^T\se\begin{bmatrix}\bar a(s)^T\\A_{12}(s)^T\end{bmatrix}$ and
$B(s)^T\se\begin{bmatrix}\bar b(s)^T\\A_{12}(s)^T\end{bmatrix}$.
Therefore, 
$ A(s)\se\begin{bmatrix}\bar a(s)&A_{12}(s)\end{bmatrix}$ and
$B(s)\se\begin{bmatrix}\bar b(s)&A_{12}(s)\end{bmatrix}$.
 \end{itemize}

\medskip

\noindent
\underline{Case $\bc=\bd$}. 
The conclusion follows applying the previous result of the case $\bu=\bv$ to the pencils $A(s)^T$ and $B(s)^T$.

\medskip

\noindent
\underline{Case $\bc\neq \bd$, $\bu\neq \bv$}. 
Assume that 
there exist homogeneous polynomials 
$\pi_1^1(s,t)\mid \dots  \mid \pi_{n}^1(s,t)$ satisfying (\ref{boca}).

\begin{enumerate}
\item[(a)]
If (\ref{cdsr}) and (\ref{sumpi1}) hold, then $q-n_1=q-n_2+1$, i.e.,
$n_1=n_2-1$, hence $n_1=n$ and $n_2=n+1$. 
From  (\ref{sumpi1}), there exists a pencil 
$A_{21}(s)\in \FF^{(p-1)\times q}$ of $\rank (A_{21}(s))=n$, homogeneous invariant factors $\pi^1_i(s, t)\mid \dots \mid \pi^1_{n}(s, t)$, column minimal indices 
$c_1\geq \dots \geq c_{q-n}$ and row  minimal indices
$v_1\geq  \dots \geq v_{p-n-1}$.
Moreover, because of (\ref{boca}), 
$$\begin{array}{rl}
\sum_{i=1}^{p-n-1}v_i&=n-\sum_{i=1}^{q-n}c_i-\sum_{i=1}^{n}\deg(\pi_i^1(s, t))\\&\leq n-\sum_{i=1}^{q-n}c_i-\sum_{i=1}^{n}\deg( \phi_i(s, t))=
\sum_{i=1}^{p-n}u_i.\end{array}$$
From  Lemma \ref{lemmacard}, we obtain
$\#\{i \; : \;\ v_i>0\}\leq \#\{i \; : \; u_i>0\}.$
Applying Lemmas \ref{lemmaDox0} and \ref{lemmaDox1}, there exist pencils
$ a(s),  b(s)\in\FF[s]^{1\times q}$ such that 
$ A(s)\se\begin{bmatrix} a(s)\\A_{21}(s)\end{bmatrix}$ and
$B(s)\se\begin{bmatrix} b(s)\\A_{21}(s)\end{bmatrix}$.
The sufficiency follows from Lemma  \ref{lemmaeq}.

\end{enumerate}

The cases (b), (c) and (d) are similar.

 \hfill $\Box$ 

\medskip

If $\FF$ is algebraically closed, the conditions of the case $\bc\neq \bd$, $\bu\neq \bv$ can be written in terms of inequalities, as stated in the next lemma. The proof is inspired by that of  \cite[Corollary 4.3]{Za97}.

\begin{lemma}\label{algclos}
  Let $\Omega_1(s, t),  \dots, \Omega_{n}(s, t),
  \Psi_1(s, t),  \dots, \Psi_{n+1}(s, t)\in \FF[s, t]$ be homogeneous polynomials such that
  $\Omega_1(s, t)\mid  \dots\mid \Omega_{n}(s, t)$, 
  $\Psi_1(s, t)\mid  \dots\mid \Psi_{n+1}(s, t)$,
  and
  \begin{equation}\label{eqintfihr1T}
\Psi_{i-1}(s, t)\mid\Omega_i(s, t)\mid\Psi_{i+1}(s, t), \quad 1\leq i \leq n.
\end{equation}
Let $x$  be a nonnegative integer.
  
If $\FF$ is an algebraically closed field, then there exist homogeneous polynomials $\pi_1^1(s,t)\mid \dots  \mid \pi_n^1(s,t)$ satisfying 
\begin{equation}\label{bocaT}
\lcm(\Omega_i(s,t), \Psi_i(s,t))\mid \pi_i^1(s,t)\mid \gcd(\Omega_{i+1}(s,t), \Psi_{i+1}(s,t)), \quad 1\leq i \leq n,
\end{equation}
and
\begin{equation}\label{sumpi1x}
\sum_{i=1}^n \deg(\pi_i^1(s,t))=x
\end{equation}
if and only if
\begin{equation}\label{bocasumx}
\sum_{i=1}^n\deg(\lcm(\Omega_i(s,t), \Psi_i(s,t)))\leq  x\leq \sum_{i=1}^n\deg(\gcd(\Omega_{i+1}(s,t), \Psi_{i+1}(s,t))).
\end{equation}
\end{lemma}

{\it Proof.}
  From (\ref{bocaT}) and (\ref{sumpi1x}), clearly we deduce (\ref{bocasumx}).

Conversely,  assume that (\ref{bocasumx}) holds.
Condition  (\ref{eqintfihr1T}) implies that
$$\lcm(\Omega_i(s,t), \Psi_i(s,t))\mid \gcd(\Omega_{i+1}(s,t), \Psi_{i+1}(s,t)), \quad 1\leq i \leq n,$$ 
hence,
$\Delta_i(s,t)=\frac{\gcd(\Omega_{i+1}(s,t), \Psi_{i+1}(s,t))}{\lcm(\Omega_i(s,t), \Psi_i(s,t))}$ 
are homogeneous polynomials.

Let $\delta_i=\deg(\lcm(\Omega_i(s,t), \Psi_i(s,t))$,  $\delta'_i=\deg(\gcd(\Omega_{i+1}(s,t), \Psi_{i+1}(s,t))), \ 1\leq i \leq n.$ From (\ref{bocasumx}) we have
$$
0\leq x-\sum_{i=1}^n\delta_i\leq\sum_{i=1}^n(\delta'_i-\delta_i).
$$
Let $z_1, \dots, z_n$ be integers such that $0\leq z_i \leq \delta'_i-\delta_i=\deg(\Delta_i(s,t))$ and  $\sum_{i=1}^nz_i= x-\sum_{i=1}^n\delta_i$.  
As $\FF$ is algebraically closed, there exists homogeneous polynomials
$\gamma_i(s,t)$ such that $\deg(\gamma_i(s,t))=z_i$ and
$\gamma_i(s,t)\mid \Delta_i(s,t)$, for $1\leq i \leq n$.

Let $\pi_i^1(s,t)= \lcm(\Omega_i(s,t), \Psi_i(s,t))\gamma_i(s, t)$, $1\leq i \leq n$. Then, $\pi^1_{i}(s,t)\mid \pi^1_{i+1}(s,t)$ for $1\leq i \leq n-1$, and they  satisfy (\ref{bocaT}) and (\ref{sumpi1x}). \hfill $\Box$ 

\begin{example}
  Let $\FF= \CC$, $n=6$, $x=3$.
  $\Omega_1(s,t)=\dots=\Omega_5(s,t)=1,\Omega_6(s,t)=s^2+t^2 $,  $\Psi_1(s,t)=\dots\Psi_5(s,t)=1,\Psi_6(s,t)=\Psi_7(s,t)=s^2+t^2 $.
  Then
  $$\lcm(\Omega_i(s,t), \Psi_i(s,t))=1, \, 1\leq i \leq 5;\;
  \lcm(\Omega_6(s,t), \Psi_6(s,t))=s^2+t^2,$$
  $$ \gcd(\Omega_{i+1}(s,t), \Psi_{i+1}(s,t))=1,  \, 1\leq i \leq 4;\;
  \gcd(\Omega_{i+1}(s,t), \Psi_{i+1}(s,t))=s^2+t^2,  \, 5\leq i \leq 6.$$
  and  (\ref{bocasumx}) holds.
  The homogeneous polynomials
$$
  \pi^1_1(s,t)=\dots =\pi^1_4(s,t)=1,\quad  \pi^1_5(s,t)\mid s+it,
  \quad \pi^1_6(s,t)= s^2+t^2
  $$
  satisfy (\ref{bocaT}) and
(\ref{sumpi1x}).
  
\end{example}

\begin{corollary}
Under the  conditions of Theorem \ref{maintheogen}, 
if $\FF$ is an algebraically closed field and
$\bc\neq \bd$, $\bu\neq \bv$,
there exists a matrix pencil $P(s)\in \FF[s]^{p \times q}$ of $\rank (P(s))=1$ such that 
$A(s)+P(s)\se B(s)$ if and only if  (\ref{eqintfihr1T}),
and one of the four following conditions hold
\begin{enumerate}
\item[(a)] (\ref{cdsr}) and
\begin{equation}\label{bocacdsr}
  \sum_{i=1}^n\deg(\lcm(\phi_i(s,t), \psi_i(s,t)))\leq  
  n-\sum_{i=1}^{q-n_1}c_i-\sum_{i=1}^{p-n_2}v_i
  \leq \sum_{i=1}^n\deg(\gcd(\phi_{i+1}(s,t), \psi_{i+1}(s,t))).
\end{equation}

\item[(b)]
(\ref{dcrs}) and
\begin{equation}\label{bocacdrs}
  \sum_{i=1}^n\deg(\lcm(\phi_i(s,t), \psi_i(s,t)))\leq  
  n-\sum_{i=1}^{q-n_2}d_i-\sum_{i=1}^{p-n_1}u_i
  \leq \sum_{i=1}^n\deg(\gcd(\phi_{i+1}(s,t), \psi_{i+1}(s,t))).
\end{equation}

\item[(c)]
  (\ref{cdsr}) y (\ref{bocacdrs}).
 \item[(d)]
  (\ref{dcrs}) y (\ref{bocacdsr}).

\end{enumerate}
  \end{corollary}

\section{Conclusions}
\label{secconclusions}

Given a  matrix pencil, regular or singular, we have completely characterized the Kronecker structure  of a  pencil obtained from  it by a perturbation of  rank one.
The  result holds over  arbitrary fields.

\bibliographystyle{acm}
\bibliography{referencesreg}

\begin{thebibliography}{10}

\bibitem{BaRo19}
{\sc {Baraga{\~n}a}, I., and {Roca}, A.}
\newblock {Fixed rank perturbations of regular matrix pencils}.
\newblock {\em arXiv e-prints\/} (Jul 2019), arXiv:1907.10657.

\bibitem{BaRo18}
{\sc {Baraga\~na}, I., and Roca, A.}
\newblock {Weierstrass structure and eigenvalue placement of regular matrix
  pencils under low rank perturbation}.
\newblock {\em SIAM Journal on Matrix Analysis and Applications 40}, 2 (2019),
  440--453.

\bibitem{Batzke14}
{\sc Batzke, L.}
\newblock Generic rank-one perturbations of structured regular matrix pencils.
\newblock {\em Linear Algebra Appl 458\/} (2014), 638--670.

\bibitem{BaMeRaRo16}
{\sc Batzke, L., Mehl, C., Ran, A., and Rodman, L.}
\newblock {Generik rank-$k$ Perturbations of Structured Matrices}.
\newblock {\em Operator Theory 255\/} (2016), 27--48.

\bibitem{TeDo07}
{\sc {De Ter\'an}, F., and Dopico, F.}
\newblock {Low rank perturbation of Kronecker structures without full rank}.
\newblock {\em SIAM Journal on Matrix Analysis and Applications 29}, 2 (2007),
  496--529.

\bibitem{TeDo16}
{\sc {De Ter\'an}, F., and Dopico, F.}
\newblock Generic change of the partial multiplicities of regular matrix
  pencils under low-rank perturbations.
\newblock {\em SIAM Journal on Matrix Analysis and Applications 37}, 3 (2016),
  823--835.

\bibitem{TeDoMo08}
{\sc {De Ter\'an}, F., Dopico, F., and Moro, J.}
\newblock {Low rank perturbation of Weierstrass structure}.
\newblock {\em SIAM Journal on Matrix Analysis and Applications 30}, 2 (2008),
  538--547.

\bibitem{Do13}
{\sc Dodig, M.}
\newblock {Completion up to a matrix pencil with column minimal indices as the
  only nontrivial Kronecker invariants}.
\newblock {\em Linear Algebra and its Applications 438\/} (2013), 3155--3173.

\bibitem{DoStEJC10}
{\sc Dodig, M., and Sto$\check{\mbox{s}}$i\'c, M.}
\newblock {On convexity of polynomial paths and generalized majorizations}.
\newblock {\em Electronic Journal of Combinatorics 17}, 1 (2010), 61.

\bibitem{DoSt13}
{\sc Dodig, M., and Sto$\check{\mbox{s}}$i\'c, M.}
\newblock On properties of the generalized majorization.
\newblock {\em Electronic Journal of Linear Algebra 26\/} (2013), 471--509.

\bibitem{DoSt14}
{\sc Dodig, M., and Sto$\check{\mbox{s}}$i\'c, M.}
\newblock The rank distance problem for pairs of matrices and a completion of
  quasi-regular matrix pencils.
\newblock {\em Linear Algebra and its Applications 457\/} (2014), 313--347.

\bibitem{DoSt19}
{\sc Dodig, M., and Sto$\check{\mbox{s}}$i\'c, M.}
\newblock The general matrix completion problem: a minimal case.
\newblock {\em SIAM Journal on Matrix Analysis and Applications 40}, 1 (2019),
  347--369.

\bibitem{Friedland80}
{\sc Friedland, S.}
\newblock {\em {Matrices: algebra, analysis and applications}}.
\newblock World Scientific, Singapore, 2016.

\bibitem{Ga74}
{\sc Gantmacher, F.}
\newblock {\em {Matrix Theory, Vols I, II}}.
\newblock Chelsea, New York, 1974.

\bibitem{GeTr17}
{\sc Gernandt, H., and Trunk, C.}
\newblock Eigenvalue placement for regular matrix pencils with rank one
  perturbations.
\newblock {\em SIAM Journal on Matrix Analysis and Applications 38}, 1 (2017),
  134--154.

\bibitem{LeMaPhTrWi18}
{\sc Leben, L., Mart\'inez-Per\'ia, F., Philipp, F., Trunk, C., and Winkler,
  H.}
\newblock Finite rank perturbations of linear relations and singular matrix
  pencils, 2018.

\bibitem{MeMeRaRo11}
{\sc Mehl, C., Mehrmann, V., Ran, A., and Rodman, L.}
\newblock Eigenvalue perturbation theory of classes of structured matrices
  under generic structured rank one perturbations.
\newblock {\em Linear Algebra and its Applications 435\/} (2011), 687--716.

\bibitem{MoDo03}
{\sc Moro, J., and Dopico, F.}
\newblock {Low rank perturbation of Jordan structure}.
\newblock {\em SIAM Journal on Matrix Analysis and Applications 25}, 2 (2003),
  495--506.

\bibitem{Ro03}
{\sc Roca, A.}
\newblock {\em Asignaci\'on de Invariantes en Sistemas de Control}.
\newblock PhD thesis, Universitat Polit\`ecnica Val\`encia, 2003.

\bibitem{Sa02}
{\sc Savchenko, S.~V.}
\newblock Typical changes in spectral properties under perturbations by a
  rank-one operator.
\newblock {\em Mathematical Notes 74}, 4 (2003), 557--568.

\bibitem{Sa04}
{\sc Savchenko, S.~V.}
\newblock On the change in the spectral properties of a matrix under
  perturbations of sufficiently low rank.
\newblock {\em Functional Analysis and Its Applications 38}, 1 (2004), 69--71.

\bibitem{Silva88_1}
{\sc Silva, F.}
\newblock {The Rank of the Difference of Matrices with Prescribed Similarity
  Classes}.
\newblock {\em Linear and Multilinear Algebra 24\/} (1988), 51--58.

\bibitem{Th80}
{\sc Thompson, R.}
\newblock Invariant factors under rank one perturbations.
\newblock {\em Canad. J. Math 32\/} (1980), 240--245.

\bibitem{Za91}
{\sc Zaballa, I.}
\newblock Pole assignment and additive perturbations of fixed rank.
\newblock {\em SIAM Journal on Matrix Analysis and Applications 12}, 1 (1991),
  16--23.

\bibitem{Za97}
{\sc Zaballa, I.}
\newblock Controllability and hermite indices of matrix pairs.
\newblock {\em International Journal of Control 68}, 1 (1997), 61--68.

\end{thebibliography}

\end{document}